\documentclass{article}

\usepackage{amsmath, amssymb}
\usepackage{graphicx, wrapfig} 
\usepackage{amsthm}
\usepackage{url}

\theoremstyle{plain}
\newtheorem{thm}{Theorem}
\newtheorem{lem}{Lemma}

\theoremstyle{definition}

\newtheorem{ex}{Example}
\newtheorem{fact}{Fact}
\newtheorem{prop}{Proposition}

\title{Boolean algebra on region crossing change and an inequality of the region unknotting number}
\author{Dawan Chumpungam\thanks{Department of General Education, KOSEN-KMITL, King Mongkut's Instituute of Technology, Ladkrabang, 1 Chalong Krung 1 Alley, Ladkrabang, Krung Thep Maha Nakhon, 10520, Thailand. Email: dawan.ch@kmitl.ac.th} \and Ayaka Shimizu\thanks{Osaka Central Advanced Mathematical Institute, 3-3-138, Sugimoto, Osaka, Japan. Email: shimizu1984@gmail.com}}
\date{\today}
\begin{document}
\maketitle

\begin{abstract}
Using Boolean algebra, we discuss the region unknotting number of a knot, and show that the region unknotting number is less than or equal to $(c+1)/2$ for any knot with crossing number $c$. 
This is a progress from $(c+2)/2$.
\end{abstract}

\section{Introduction}

\begin{figure}[h]
	\centering
	\includegraphics[width=35mm]{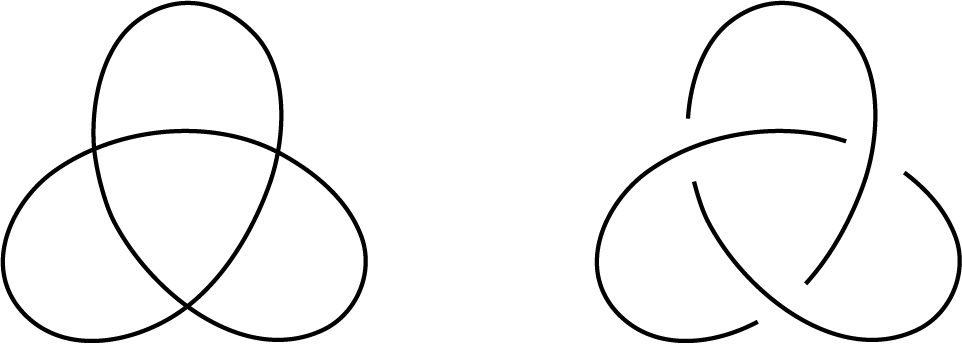}
	\caption{A knot projection and knot diagram.}
	\label{fig-kp}
\end{figure}

{\it A knot} is an embedding of a circle in $S^3$. 
{\it A knot projection} is a projection of a knot on $S^2$ such that any intersection is a double point where two arcs intersect transversely. 
{\it A knot diagram} is a projection of a knot with ``crossing information'' to represent a knot, as shown in Figure \ref{fig-kp}. 
Each knot projection with $c$ crossings divides $S^2$ into $c+2$ connected parts\footnote{We can see that the number of regions of a knot projection is larger by two than the number of crossings, by considering Euler's characteristic. See, for example, \cite{ahara-suzuki}.} and we call each part a {\it region}.

The {\it crossing number of a knot diagram $D$}, denoted by $c(D)$, is the number of crossings. 
The {\it crossing number of a knot $K$}, denoted by $mc(K)$, is the minimal value of $c(D)$ for all diagrams $D$ of $K$. 
A {\it trivial knot} is the knot with $mc(K)=0$. 
A knot diagram $D$ of $K$ is said to be {\it minimal} when $c(D) = mc(K)$.

\begin{figure}[h]
	\centering
	\includegraphics[width=40mm]{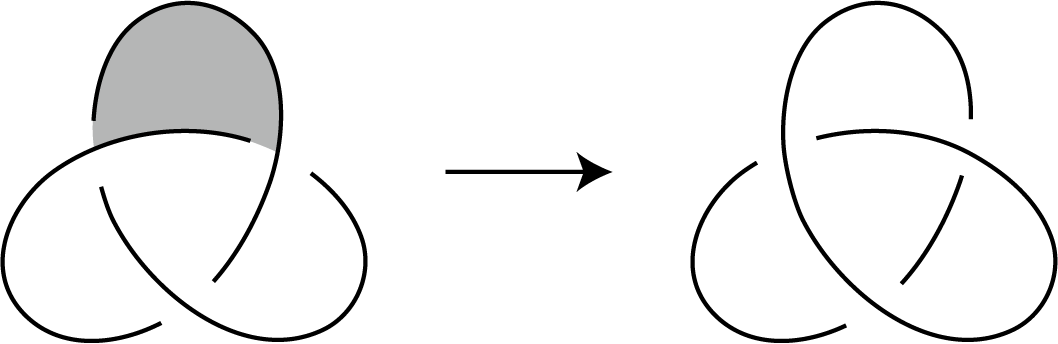}
	\caption{A region crossing change at the shaded region.}
	\label{fig-rcc}
\end{figure}

A {\it crossing change} at a crossing of $D$ is a local move changing the crossing information of the crossing. 
A {\it region crossing change}\footnote{RCC was defined by Kengo Kishimoto in a seminar at Osaka City University in 2010.}, RCC, on a region of $D$ is a set of crossing changes at all the crossings on the boundary of the region. 
See Figure \ref{fig-rcc}.

Local moves, such as a crossing change, region crossing change, etc. can be used to measure a complexity or a knot of a knot diagram. 
The {\it region unknotting number} of a knot diagram $D$, denoted by $u_R(D)$, is the minimum number of RCCs on regions which are needed to transform $D$ into a diagram of the trivial knot. 
The {\it region unknotting number}\footnote{It is shown in \cite{aida} that any knot can be unknotted by a single ``$n$-gon move'', and this means that any nontrivial knot has a diagram $D$ with $u_R(D)=1$. The point of $mu_R(K)$ is that we have to choose a {\it minimal diagram} and we cannot retake the diagram during RCCs.} of a knot $K$, denoted by $mu_R(K)$, is the minimal value of $u_R(D)$ for all minimal diagrams $D$ of $K$. 
The following theorem was shown in \cite{rcc-uo}.

\begin{thm}
	$mu_R(K) \leq \dfrac{c+2}{2}$ holds for any knot $K$ with crossing number $c$.
	\label{thm-ru2}
\end{thm}

\noindent Some sharper inequalities have been given for some specific types of knots (see, for example, \cite{rcc-un, rcc-uo, rcc-torus, rcc-2bridge}). 
In this paper, we show the following inequality.

\begin{thm}
	$mu_R(K) \leq \dfrac{c+1}{2}$ holds for any knot $K$ with crossing number $c$. 
	\label{thm-ru1}
\end{thm}

\noindent For each fixed knot diagram, the region unknotting number can be determined in the following way. 
All the sets of regions which unknot the knot diagram by RCC can be found by solving a system of linear equations (see, for example, \cite{ahara-suzuki}), and by looking at the minimum value of the cardinalities of the sets, we can find the region unknotting number of the diagram. 
However, no algorithms to find the number directly have been found at the moment. 
As a side note, various kinds of problems regarding the minimum number have been studied as an optimization problem, and various problems and methods using lattices are actively explored (see, for example, \cite{opt-1, opt-2, opt-3}). 
In this paper, the authors attempt to describe RCC in terms of the Boolean algebras which are suitable with RCC and would be useful to address the optimization problem of the region unknotting number in the future studies.

\noindent The rest of the paper is organized as follows. 
In Section \ref{section-pr}, some properties of RCC which are required in this paper are listed. 
In Section \ref{section-b}, the definition and useful facts on Boolean algebra are listed. 
In Section \ref{section-rcc-b}, RCC and its properties are described in terms of Boolean algebra. 
In Section \ref{section-rcc}, the proof of Theorem \ref{thm-ru1} is given.

\section{Properties of region crossing change}
\label{section-pr}

In this section, some properties of the region crossing change are listed. 
First, the result on the RCCs does not depend on the order of regions; 
Around a crossing, if an odd (resp. even) number of regions are applied RCC, then the crossing changes (resp. does not change) as a result. 
Therefore, we can consider RCCs as an operation on a set of regions. Let $S, \ T$ be sets of regions of a knot diagram. 
Applying RCCs on $S$ and then $T$ is equivalent to the RCC on the symmetric difference $S \oplus T =(S \cup T) \setminus (S \cap T)$ because RCCs twice on the same region are cancelled as a result. 
Moreover, let $C_1$ (resp. $C_2$) be the set of crossings which are affected by RCC on a set of regions $R_1$ (resp. $R_2$). 
Then, the RCC on $R_1 \oplus R_2$ affects $C_1 \oplus C_2$ because crossing-changes twice are cancelled, too.

\begin{figure}[h]
	\centering
	\includegraphics[width=50mm]{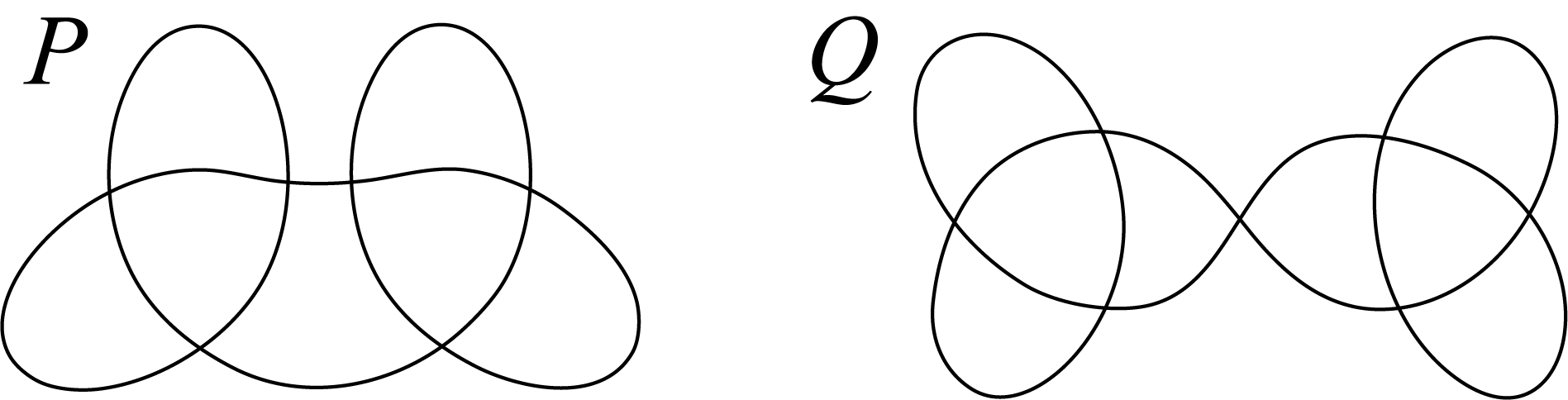}
	\caption{An irreducible projection $P$ and reducible projection $Q$.}
	\label{fig-red}
\end{figure}

A knot projection $P$ is said to be {\it irreducible} if any crossing of $P$ has four disjoint regions around the crossing. 
On the other hand, a knot projection $Q$ is said to be {\it reducible} if $Q$ has a crossing which has only three regions around the crossing. See Figure \ref{fig-red}.

\begin{figure}[h]
	\centering
	\includegraphics[width=40mm]{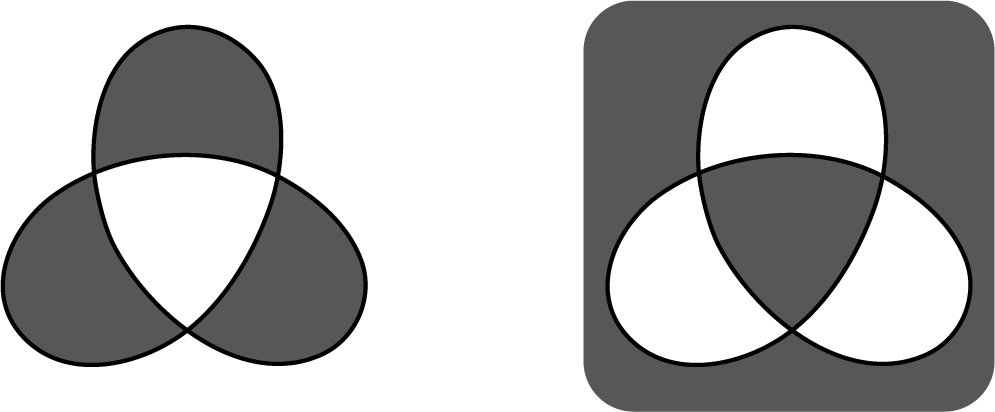}
	\caption{Checkerboard coloring to a knot projection.}
	\label{fig-check}
\end{figure}

A {\it checkerboard coloring} is a coloring of two colors, say black and white, to each region of a knot projection such that each pair of regions sharing an edge have different color (see Figure \ref{fig-check}). 
We note that there are two choices of a checkerboard coloring for each knot projection.

Let $P$ be an irreducible knot projection with a checkerboard coloring. 
Let $B$ (resp. $W$) be the set of all regions of the black (resp. white) regions. 
The following lemmas were shown in \cite{rcc-uo}. 

\begin{lem}
	For an irreducible knot projection, RCC on neither $B$ nor $W$ affects any crossing as a result.\footnote{These kinds of sets are called ineffective sets and studied in various settings. For example, see \cite{russell, inoue, rinno}. }
	\label{lem-bw}
\end{lem}

\begin{lem}
	On an irreducible knot projection, $B_S$ and $B \oplus B_S$ (resp. $W_S$ and $W \oplus W_S$) result in the same on RCC for any subset $B_S \subset B$ (resp. $W_S \subset W$). 
	\label{lem-bw-s}
\end{lem}

\begin{lem}
	For any set of regions $S$ of an irreducible knot projection, the four sets $S$, $S \oplus B$, $S \oplus W$ and $S \oplus (B \oplus W)$ result in the same on RCC. 
	\label{lem-bw-4}
\end{lem}

\noindent For a set of regions $S$, we call the other three sets of regions $S \oplus B$, $S \oplus W$, $S \oplus (B \oplus W)$ the {\it BW-complements of $S$}. 
The following theorem was shown in \cite{rcc-uo}.

\begin{thm}
	Any crossing change on a knot diagram can be realized by RCCs\footnote{Theorem \ref{thm-rcc-unknot} implies that a region crossing change is an {\it unknotting operation} on a knot diagram. See \cite{cheng-gao, cheng} for links. See \cite{FNS, HSS, rinno} for spatial-graphs.}. 
	\label{thm-rcc-unknot}
\end{thm}

\begin{figure}[h]
	\centering
	\includegraphics[width=90mm]{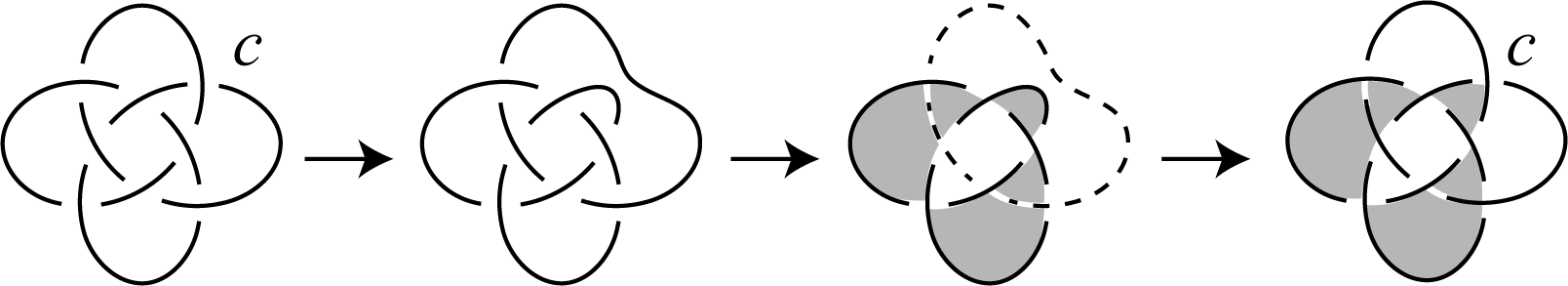}
	\caption{We can find a set of regions which changes a crossing $c$ in the following procedure; Splice $c$ to obtain a two-component link diagram. Apply a checkerboard coloring to one of the components. Then, the corresponding regions change $c$ by RCC.}
	\label{fig-alg}
\end{figure}

\noindent An algorithm to find such a set of regions for an irreducible knot projection is shown in Figure \ref{fig-alg}.

The effect of RCC can be seen in a system of linear equations; 
For a knot projection with $c$ crossings, give labels $c_1 , c_2 , \dots , c_c$ to the crossings and $R_1 , R_2, \dots , R_{c+2}$ to the regions in any order. 
The {\it region choice matrix} \cite{ahara-suzuki} (or transpose of the {\it incidence matrix} \cite{cheng-gao}) is a $c \times (c+2)$ matrix $M=(m_{ij})$ defined by

\vspace{-4mm}
\begin{align*}
	m_{i j}= \left\{
	\begin{array}{ll}
		1 & \text{when } c_i \text{ is incident to } R_j \\
		0 & \text{otherwise}
	\end{array}
	\right.
\end{align*}

\begin{figure}[h]
	\centering
	\includegraphics[width=60mm]{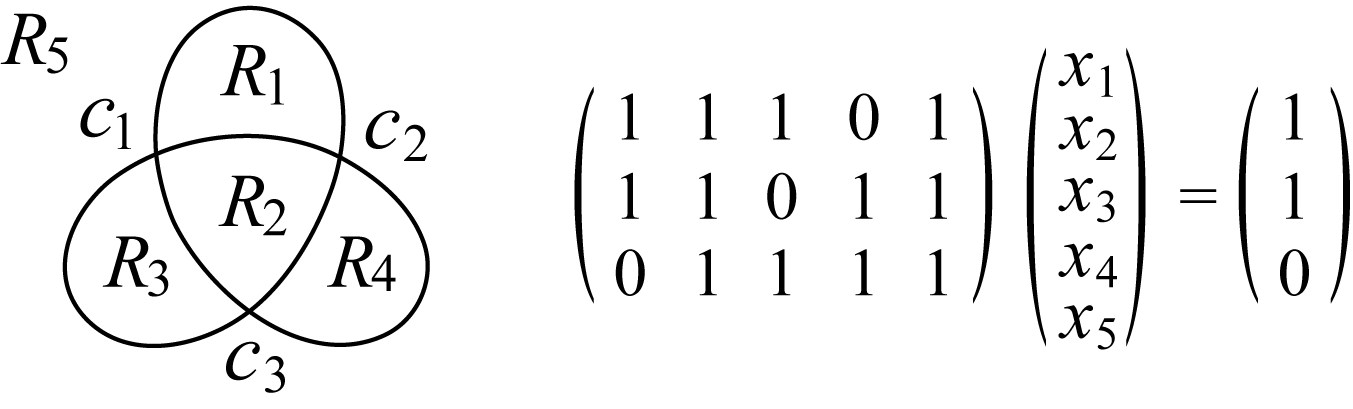}
	\caption{By solving the system of linear equations, we have the solution $(x_1, x_2, x_3, x_4, x_5)= (1,0,0,0,0), (1,1,0,0,1), (0,0,1,1,0), (0,1,1,1,1)$, namely, $\{ R_1 \} , \{ R_1 , R_2 , R_5 \} , \{ R_3 , R_4 \} ,  \{ R_2 , R_3 , R_4 , R_5 \}$ to change the crossings $c_1$ and $c_2$.}
	\label{fig-matrix}
\end{figure}

\noindent If one wants to change the $i$th crossings, prepare a vector $\boldsymbol{b} \in {\mathbb{R}}^c$ with the $i$th elements $1$ and the others $0$, and solve $M \boldsymbol{x}= \boldsymbol{b} \ ( \text{mod } 2)$. 
See Figure \ref{fig-matrix}. 
Since the $c \times (c+2)$ region choice matrix $M$ is full-rank for any knot projection (\cite{cheng-gao}) and we work on $\mathbb{Z}/2 \mathbb{Z}$, we have exactly four corresponding sets of regions for any set of crossings \cite{kawauchi}.

\noindent The following lemma was shown in \cite{rcc-un} (see \cite{k-col, fukushima, jong} for more details).

\begin{lem}
	Let $P$ be an irreducible knot projection with a checkerboard coloring. 
	For any pair of black and white regions, we can make any crossing change on $P$ by RCCs without using the two regions. 
	\label{lem-two}
\end{lem}

\begin{figure}[h]
	\centering
	\includegraphics[width=50mm]{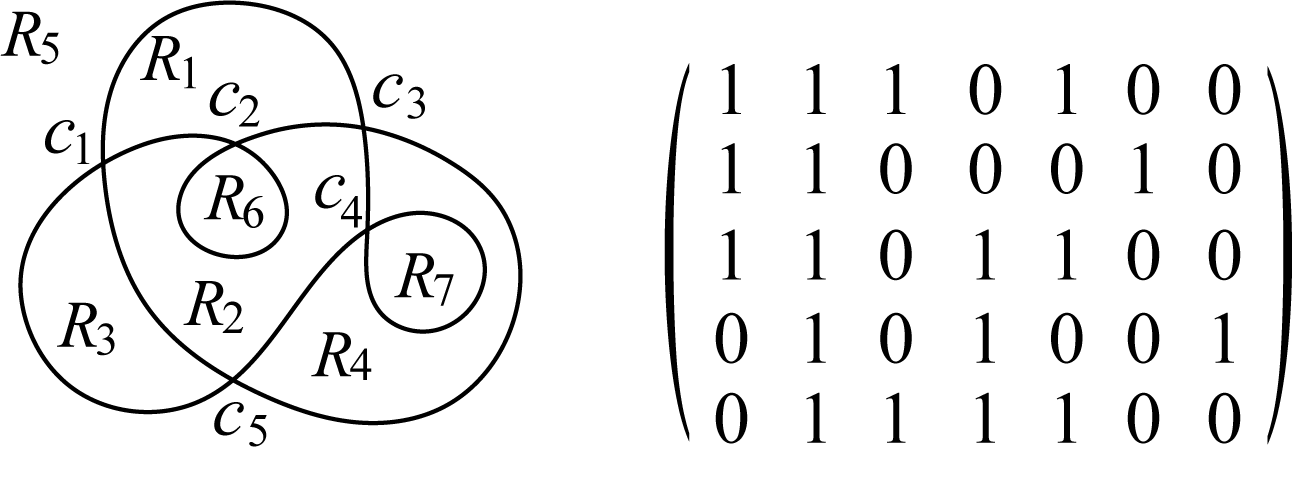}
	\caption{We cannot make a single crossing change at $c_1$ by RCC on any subset of $\{ R_1 , R_2 , R_3 , R_4 , R_5 \}$; The square matrix of the first to fifth columns in the region choice matrix is not invertible.}
	\label{fig-counter}
\end{figure}

\noindent We remark that Lemma \ref{lem-two} does not hold for reducible knot projections. 
See Figure \ref{fig-counter}.

\section{Boolean algebra}
\label{section-b}

As shown in Section \ref{section-rcc-b}, the Boolean algebra is suitable for the discussion on RCC. 
In this section, we prepare some terms of Boolean algebra based on \cite{judson}.

A poset $(L, \preceq)$ is said to be a {\it lattice} if $L$ has the least upper bound $a \vee b$ and the greatest lower bound $a \wedge b$ for any pair of elements $a, b \in L$. 
An element $I$ is said to be a {\it largest element} if $a \preceq I$ holds for any element $a \in L$. 
An element $O$ is said to be a {\it smallest element} if $O \preceq a$ holds for any element $a \in L$. 
A lattice is {\it distributive} if the distributive law $a \wedge (b \vee c)= (a \wedge b) \vee (a \wedge c)$ holds. 
A lattice $L$ with $I$ and $O$ is said to be {\it complemented} if ${}^{\exists} a' \in L$ such that $a \vee a' =I$ and $a \wedge a' = O$ for each $a \in L$. 
A {\it Boolean algebra} is a distributive and complemented lattice $B$ with a largest element $I$ and a smallest element $O$. 
The following fact is well known.

\vspace{2mm} 

\begin{fact}
	A set $B$ is a Boolean algebra if and only if $B$ has binary operations $\vee$ and $\wedge$ which satisfy the following five conditions. \\
	1. $a \vee b = b \vee a$ and $a \wedge b = b \wedge a$ for ${}^{\forall}a, b \in B$. \\
	2. $a \vee (b \vee c) = (a \vee b) \vee c$ and $a \wedge (b \wedge c)= (a \wedge b) \wedge c$ for ${}^{\forall}a, b, c \in B$. \\
	3. $a \wedge (b \vee c) = (a \wedge b) \vee (a \wedge c)$ and $a \vee (b \wedge c)= (a \vee b) \wedge (a \vee c)$ for ${}^{\forall}a, b, c \in B$. \\
	4. ${}^{\exists} I$ and ${}^{\exists}O$ such that $a \vee O=a$ and $a \wedge I=a$ for ${}^{\forall}a \in B$. \\
	5. For ${}^{\forall}a \in B$, ${}^{\exists}a' \in B$ such that $a \vee a' = I$ and $a \wedge a' = O$. 
\end{fact}

\noindent The following fact is also well known.

\begin{fact}
	On a Boolean algebra, $ a \preceq b \Leftrightarrow a \wedge b' = O \Leftrightarrow a' \vee b = I$.
	\label{fact-b}
\end{fact}

\noindent The {\it power set} $(P(S), \cup, \cap, \subseteq)$ of a set $S$ is the set of all subsets of $S$, and is a typical example of a lattice or a Boolean algebra.

\begin{figure}[h]
	\centering
	\includegraphics[width=70mm]{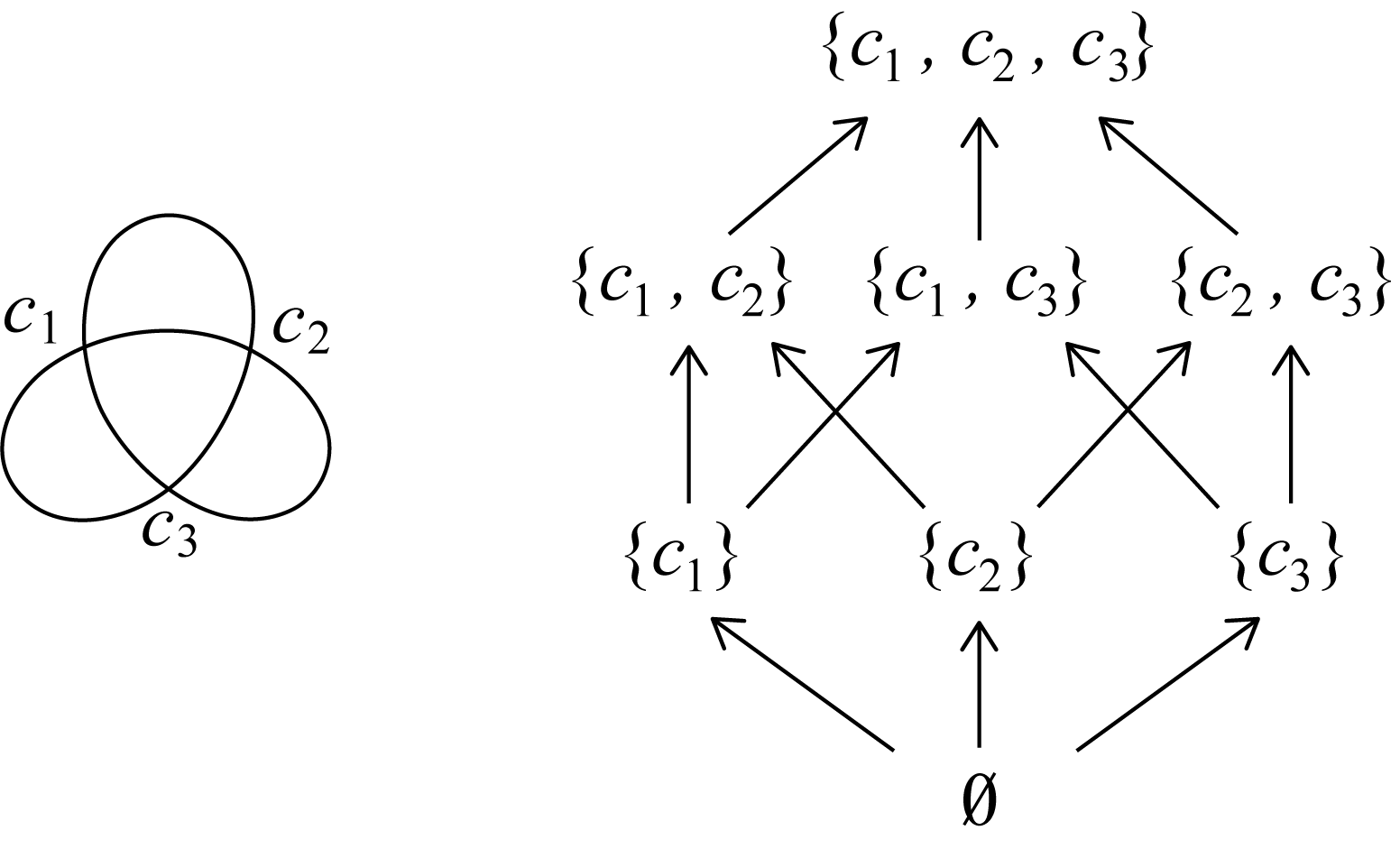}
	\caption{A Boolean algebra $P(C)$ with $a \vee b = a \cup b$ and $a \wedge b = a \cap b$.}
	\label{fig-pc}
\end{figure}

\begin{ex}
	Let $C$ be the set of all crossings of a knot projection. 
	Then the power set $P(C)$ is a Boolean algebra with the binary operations $a \vee b = a \cup b$ and $a \wedge b = a \cap b$, the union and intersection of sets. 
	Here, the complement $a'$ of $a$ is the complement $C \setminus a$. (See Figure \ref{fig-pc}.)
	\label{ex-pc}
\end{ex}

\section{Region crossing change in terms of Boolean algebra}
\label{section-rcc-b}

In this section, we describe the region crossing change using Boolean algebras.

Let $C$ (resp. $R$) be the sets of all crossings (resp. regions) of a knot projection. 
As mentioned in Example \ref{ex-pc}, the power set $P(C)$, and similarly $P(R)$ are Boolean algebra with $a \vee b = a \cup b$, $a \wedge b = a \cap b$, and $a'=C \setminus a$ for $P(C)$ and $a'= R \setminus a$ for $P(R)$. 
We note that the order of $P(C)$ is $2^c$ and that of $P(R)$ is $2^{c+2}$ for each knot projection of $c$ crossings.
We denote the largest element of $P(C)$ (resp. $P(R)$) by $I^c$ (resp. $I^R$), and the smallest element of $P(C)$ (resp. $P(R)$) by $O^c$ (resp. $O^R$). 
Note that $I^c =C$, $I^R =R$, $O^c = O^R = \emptyset$.

The effect of RCC can be considered as a Boolean function $\phi : P(R) \to P(C)$. 
From the region choice matrix, we see that we have exactly four elements in $P(R)$ which has the image $a$ for each $a \in P(C)$.

Give a checkerboard coloring to a knot projection. 
Let $B$ (resp. $W$) be the set of all black (resp. white) regions. 
Then $B$, $W$ are elements of $P(R)$ and Lemma \ref{lem-bw} is rephrased as follows.

\begin{lem}
	$\phi(B)= \phi(W)=O^c$ for any irreducible knot projection.
\end{lem}

Now we introduce one more useful operation $a \oplus b = ( a \wedge b') \vee (a' \wedge b)$. 
Note that $a \oplus b$ is the symmetric difference $(a \cup b) \setminus (a \cap b)$ of sets here, and $\phi(a \oplus b)= \phi(a) \oplus \phi(b)$ holds for ${}^{\forall}a, b \in P(R)$. 
For the BW-complements, Lemma \ref{lem-bw-4} is rephrased as follows:

\begin{lem}
	For any $S \in P(R)$, $\phi(S)=\phi(S \oplus B)=\phi(S \oplus W)=\phi(S \oplus I^R)$ on an irreducible knot projection. 
	\label{lem-bwc}
\end{lem}

\section{Upper bound for the region unknotting number}
\label{section-rcc}

In this section, we prove Theorem \ref{thm-ru1}.

\subsection{Monotone diagram}

\begin{figure}[h]
	\centering
	\includegraphics[width=55mm]{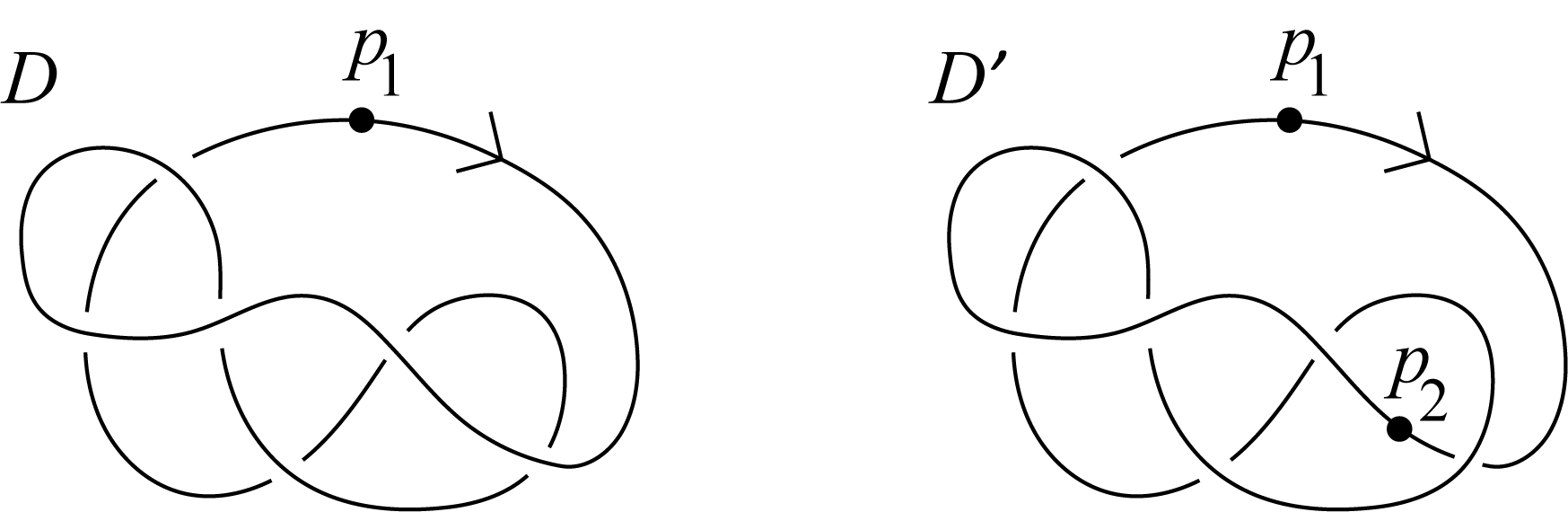}
	\caption{Monotone diagrams $D$ from $p_1$ and $D'$ from $p_2$.}
	\label{fig-mon}
\end{figure}

Let $D$ be an oriented knot diagram. 
Take a base point $p$ on $D$, avoiding crossings. 
$D$ is said to be {\it monotone from $p$} if we meet any crossing as an over-crossing first when we travel $D$ from $p$ in the given orientation \cite{lecture}. 
See Figure \ref{fig-mon}. 
It is known that any monotone knot diagram represents the trivial knot. The following lemma will be used in the proof of Theorem \ref{thm-ru1}.

\begin{lem}
	Let $D$ be a monotone diagram from a base point $p_1$. 
	Let $D'$ be the diagram obtained from $D$ by a crossing change at the crossing just after $p_1$ with the orientation. 
	Then, $D'$ is monotone from a base point $p_2$ just after the crossing. 
	\label{lem-monotone}
\end{lem}

\subsection{Equilibrium set of regions}

\noindent From this subsection, we work on the Boolean algebra again introduced in Section \ref{section-rcc-b}. 
Let $B$ (resp. $W$) $\in P(R)$ be the set of all black (resp. white) regions for a knot projection with a checkerboard coloring. 
A set of regions $S \in P(R)$ is said to be {\it equilibrium} if $|S \wedge B|= |B|/2$ and $|S \wedge W|=|W|/2$, where $|a|$ denotes the order of (or the number of regions in) $a$. 
For example, the set $S$ in Figure \ref{fig-p0-12} is equilibrium. 
Note that the equilibriumness does not depend on the choice of a checkerboard coloring. 
By definition, we have the following.

\begin{lem}
	There exists an equilibrium element in $P(R)$ if and only if both $|B|$ and $|W|$ are even with a checkerboard coloring.
	\label{lem-eq-even}
\end{lem}

\begin{figure}[h]
	\centering
	\includegraphics[width=80mm]{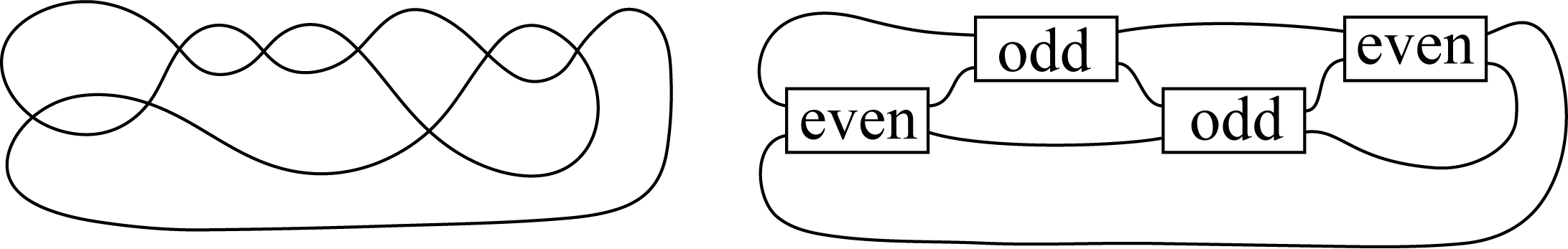}
	\caption{The knot projection on the left-hand side has even numbers of $|B|$ and $|W|$. This knot projection is known as a projection of a ``two-bridge knot'' or a ``rational knot'' $T(2,3,1,2)$ (see, for example, \cite{adams}). Generally, any rational knot projection of type $T(e,o,o,e)$, shown in the right-hand side, has even numbers of $|B|$ and $|W|$, where $e$ is an even number and $o$ is an odd number. Besides, any rational knot projection of type $T(e,o,o,o,o)$, $T(o,o,o,o,o,o)$, $T(o,e,o,o,o,e)$, $T(e,e,e,o,o,e)$ or $T(e,o,e,e,o,e)$ has even numbers of $|B|$ and $|W|$, too.}
	\label{fig-bw-even}
\end{figure}

\noindent We remark that there are infinitely many knot projections which have even numbers of $|B|$ and $|W|$. 
See Figure \ref{fig-bw-even}. By definition, we have the following.

\begin{lem}
	If $S$ is equilibrium, then the BW-compliments $S \oplus B$, $S \oplus W$, $S \oplus I^R$ are also equilibrium. 
	\label{lem-eq-comp}
\end{lem}

\subsection{Proof of Theorem \ref{thm-ru2}}

In this subsection, we show the proof of Theorem \ref{thm-ru2} given in \cite{rcc-uo}\footnote{Although the proof of Theorem \ref{thm-ru2} was improved in \cite{unknotting} to be more simple without using checkerboard coloring, we intentionally chose the original one because we will use the idea of checkerboard coloring in the proof of Theorem \ref{thm-ru1}.}, in terms of Boolean algebra.

\vspace{3mm}
\noindent Proof of Theorem \ref{thm-ru2}. 

\begin{figure}[h]
	\centering
	\includegraphics[width=90mm]{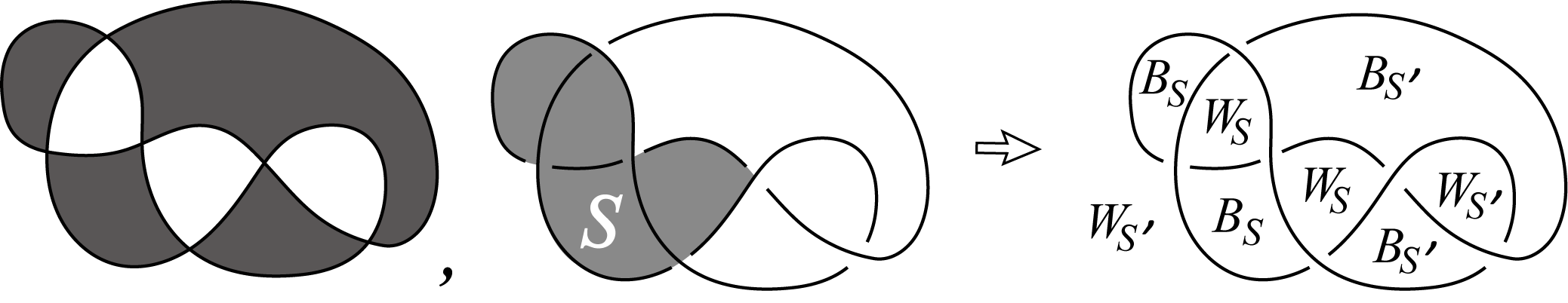}
	\caption{$B_S, B_{{S}'}, W_S$ and $W_{{S}'}$.}
	\label{fig-p0-12}
\end{figure}

Let $D$ be a minimal knot diagram of a knot $K$ with $c(D) = mc(K) = c$. 
Note that $D$ is irreducible. 
Let $R$ be the set of all the regions of $D$. 
Note that $|R|=c+2$. 
Let $S \in P(R)$ be a set of regions of $D$ which makes $D$ unknotted by RCC. 
By Lemma \ref{lem-bwc}, the BW-complements with a checkerboard coloring $S \oplus B$, $S \oplus W$ and $S \oplus I^R$ also make $D$ unknotted. 
Let $S'$ be the complement $R \setminus S$, and let $B_S = B \wedge S$, $B_{{S}'} = B \wedge S'$, $W_S = W \wedge S$, $W_{{S}'}= W \wedge S'$ (see Figure \ref{fig-p0-12}). 
Then $S= B_S \vee W_S$, $S \oplus B = B_{{S}'} \vee W_S$, $S \oplus W = B_S \vee W_{{S}'}$ and $S \oplus I^R = B_{{S}'} \vee W_{{S}'}$, and therefore $|S|= | B_S | +| W_S| $, $| S \oplus B | = | B_{{S}'} |+| W_S|$, $|S \oplus W| = | B_S  | + | W_{{S}'}| $ and $| S \oplus I^R | = |B_{{S}'} |+|  W_{{S}'}|$. 
Since $|B_S | + |B_{{S}'} | + | B_W | + |B_{{S}'}| = |B| +|W| = c+2$, we have $|S|+|S \oplus B| + |S \oplus W | + | S \oplus I^R|= 2(c+2) = 4 \times (c+2)/2$. 
Hence at least one of them is less than or equal to $(c+2)/2$. 
Therefore, $mu_R(K) \leq u_R(D) \leq (c+2)/2$. $\square$

\vspace{3mm}

\noindent Next, we show the following lemma.

\begin{lem}
	If a knot diagram $D$ with $c$ crossings has no equilibrium sets, then $u_R(D)<(c+2)/2$ holds.
	\label{lem-equil}
\end{lem}

\noindent Proof. 
Let $S \in P(R)$ be any set of regions such that RCC on $S$ (and then $S \oplus B$, $S \oplus W$, $S \oplus I^R$ by Lemma \ref{lem-bwc}) transforms $D$ into a diagram of the trivial knot. 
Let $B_S, B_{{S}'}, W_S, W_{{S}'}$ be the sets in the same setting to the proof of Theorem \ref{thm-ru2}. 
Then at least one of $|S|=|B_S |+|W_S |$, $|S \oplus B| = |B_{{S}'} | + |W_S|$, $|S \oplus W | = |B_S | + |W_{{S}'}|$ and $| S \oplus I^R | = |B_{{S}'}| + | W_{{S}'}|$ is less than $|B|/2 + |W|/2$ because $S$ and the BW-complements are not equilibrium. 
Hence, $u_R(D) < (|B|+|W|)/2 = (c+2)/2$. 
$\square$

\vspace{3mm}

\subsection{Proof of Theorem \ref{thm-ru1}}

We prove Theorem \ref{thm-ru1}. 

\vspace{3mm}

\noindent Proof of Theorem \ref{thm-ru1}. Let $D$ be a minimal diagram of a knot $K$ with $c(D) = mc(K) = c$ and $u_R(D) = mu_R(K)$. 
Note that $D$ is irreducible. 
Let $B$ (resp. $W$) $\in P(R)$ be the set of all black (resp. white) regions. 
If $|B|$ or $|W|$ is odd, we have $u_R(D) < (c+2)/2$ by Lemma \ref{lem-equil} because any set of regions cannot be equilibrium by Lemma \ref{lem-eq-even}. 
Assume $|B|=2m$ and $|W|=2n$ for some $m, n \in \mathbb{N}$.

\begin{figure}[h]
	\centering
	\includegraphics[width=25mm]{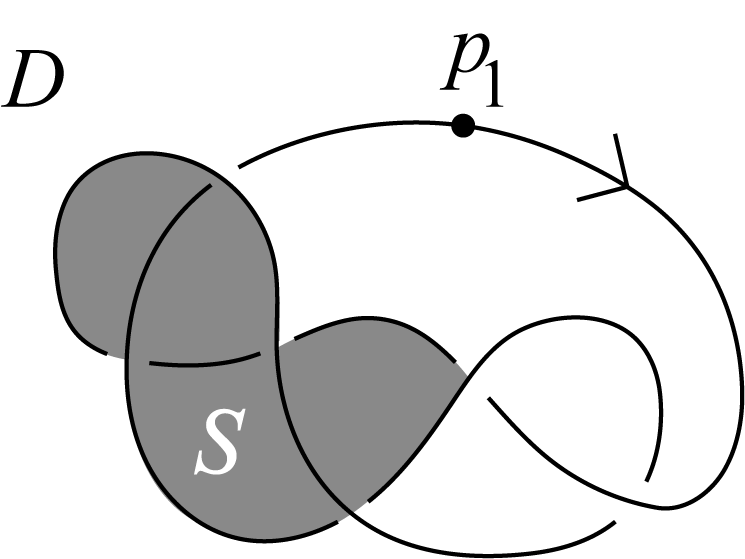}
	\caption{The set $S$ makes $D$ monotone from $p_1$ by RCCs.}
	\label{fig-p0}
\end{figure}

Procedure 0: Give an orientation to $D$ and take a base point $p_1$. 
Let $S \in P(R)$ be a set of regions which make $D$ monotone from $p_1$ by RCC. (See Figure \ref{fig-p0}.)

Let $B_S =B \wedge S$, $B_{{S}'}= B \wedge S'$, $W_S = W \wedge S$, $W_{{S}'} = W \wedge S'$, where $S'$ is the complement of $S$. (See Figure \ref{fig-p0-12}.) 
Note that the four sets $S= B_S \vee W_S$, $B_S \vee W_ {{S}'} $, $B_{{S}'} \vee W_S$, $B_{{S}'} \vee W_{{S}'}$ are BW-complements each other; they all make $D$ monotone from $p_1$ by RCC.

If $|B_S| \neq m$ or $|W_S| \neq n$, i.e., $S$ is not equilibrium, we have $u_R(D)< m+n = (c+2)/2$ by taking an appropriate BW-complement. 
If $|B_S|=m$ and $|W_S|=n$ (then also $|B_{{S}'} |=m$ and $|W_{{S}'} |=n$), i.e., $S$ is equilibrium, proceed to the next step.

\begin{figure}[h]
	\centering
	\includegraphics[width=25mm]{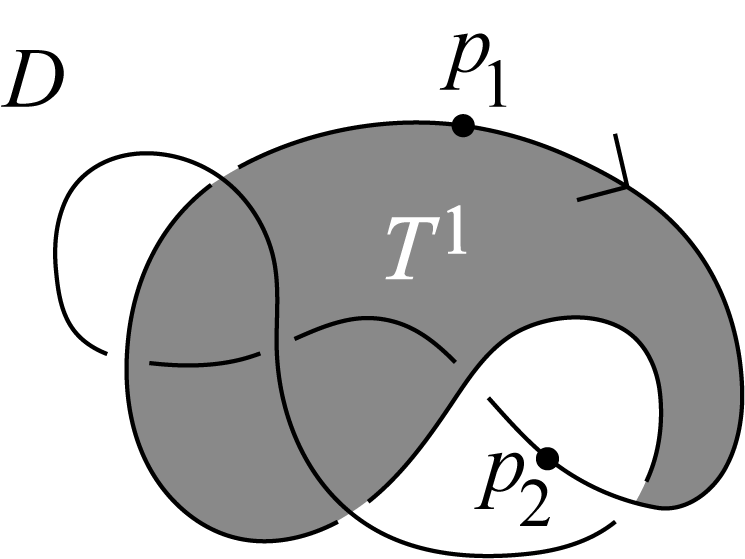}
	\caption{$p_2$ and $T^1$.}
	\label{fig-p1}
\end{figure}

Procedure 1: Take another base point $p_2$ just after a crossing from $p_1$ with the orientation. 
Let $T^1 \in P(R)$ be a set of regions such that $T^1$ changes the single crossing between $p_1$ and $p_2$ by RCC, namely, $\phi(T^1)$ is the crossing between $p_1$ and $p_2$. (See Figure \ref{fig-p1}.) 
Then $S_1 = S \oplus T^1$ makes $D$ monotone from $p_2$ by RCC by Lemma \ref{lem-monotone}.

If $|T^1 \wedge B_S| \neq |T^1 \wedge B_{{S}'} |$ or $|T^1 \wedge W_S| \neq |T^1 \wedge W_{{S}'} |$, we have $u_R(D)< (c+2)/2$ because $S_1$ is not equilibrium. 
If $|T^1 \wedge B_S| = |T^1 \wedge B_{{S}'} |$ and $|T^1 \wedge W_S| = |T^1 \wedge W_{{S}'} |$, then $S_1$ is still equilibrium. 
Moreover, the three BW-compliments are also equilibrium. 
Proceed to the next step.

Procedure 2: Take another base point $p_3$ just after a crossing from $p_2$. 
Let $T^2$ be an element of $P(R)$ such that $\phi(T^2)$ is the crossing between $p_2$ and $p_3$. 
Then $S_2 = S_1 \oplus T^2 = (S \oplus T^1 ) \oplus T^2$ makes $D$ monotone from $p_3$ by RCCs.

Let $B_{S_1} = B \wedge S_1$, $B_{{S_1}'} = B \wedge {S_1}'$, $W_{S_1} = W \wedge S_1$, $W_{{S_1}'} = W \wedge {S_1}'$, where ${S_1}'$ is the complement of $S_1$. 
If $|T^2 \wedge B_{S_1}| \neq |T^2 \wedge B_{{S_1}'}|$ or $|T^2 \wedge W_{S_1}| \neq |T^2 \wedge W_{{S_1}'}|$, we have $u_R(D)< (c+2)/2$ because $S_2$ is not equilibrium. 
If $|T^2 \wedge B_{S_1}| = |T^2 \wedge B_{{S_1}'} |$ and $|T^2 \wedge W_{S_1}| = |T^2 \wedge W_{{S_1}'}|$, then $S_2$ and the BW-compliments are equilibrium. 
Proceed to the next step.

Continue the procedure until finding the number $k$ such that $|T^k \wedge B_{S_{k-1}}| \neq |T^k \wedge B_{{S_{k-1}}'}|$ or $|T^k \wedge W_{S_{k-1}}| \neq |T^k \wedge W_{{S_{k-1}}'} |$. 
There definitely exists such $k < 2c$. 
Otherwise, it results that any $T \in P(R)$ (and its BW-complements) with $| \phi (T)|=1$ has even numbers of both black and white regions. 
Take $c$ elements $T_1 , T_2 , \dots , T_c \in P(R)$ with $\phi (T_1) = c_1 , \phi (T_2)=c_2 , \dots , \phi(T_c)=c_c$. 
Then for any $a= c_{i_1} \oplus c_{i_2} \oplus \dots \oplus c_{i_k} \in P(C)$, $b=T_{i_1} \oplus T_{i_2} \oplus \dots \oplus T_{i_k}$, which satisfies $\phi (b)=a$, and its BW-compliments have even numbers of black and white regions. 
We remark here that $|r \oplus s|$ is even when both $|r|$ and $|s|$ are even since $|r \oplus s |= |(r \vee s) \setminus (r \wedge s)| = |r| + |s| - |r \wedge s| - |r \wedge s|$. 
However, we have lots of counter-examples. 
For example, take a single region, which has just one black or white region.

Hence, $u_R(D)< (c+2)/2$ holds for any irreducible knot diagram, and $mu_R(K) \leq (c+1)/2$ holds for any knot $K$. $\square$

\section*{Appendix: Boolean equivalence between crossings and regions}

In this appendix, we show that the set of all crossings of a knot projection and a set of all regions except for an appropriate pair of regions\footnote{Again, see \cite{k-col, fukushima, jong} for more details.} are Boolean-algebraically equivalent.

Let $P$ be an irreducible knot projection with $c$ crossings. 
Give a checkerboard coloring to $P$. 
Let $b$ (resp. $w$) be a region colored black (resp. white). 
Let $S$ be the set of all the regions of $P$ except $b$ and $w$. 
As for the map $\phi$ introduced in Section \ref{section-rcc-b}, we have the following proposition.

\begin{prop}
	The map $\phi: P(S) \to P(C)$ is bijecive. 
	\label{prop-bi}
\end{prop}

\noindent Proof. 
The map $\phi: P(S) \to P(C)$ is surjective by Lemma \ref{lem-two}. 
Then, $\phi: P(S) \to P(C)$ is injective since $|P(S)|=|P(C)|$.
$\square$

\vspace{3mm}

\noindent Now we construct a Boolean algebra for $P(S)$.

\begin{prop}
	The power set $P(S)$ is a Boolean algebra with the binary operations $a {\vee}^r b = \phi^{-1} (\phi(a) \vee \phi(b))$ and $a {\wedge}^r b = \phi^{-1} (\phi(a) \wedge \phi(b))$. 
\end{prop}

\noindent Proof. 
It is sufficient to check the five conditions of Fact \ref{fact-b}. For ${}^{\forall}a, b, c \in P(S)$, \\

\noindent 1. $a {\vee}^r b = b {\vee}^r a$ and $a {\wedge}^r b = b {\wedge}^r a$ hold by definition. \\

\noindent 2.
$\begin{aligned}[t]
	a {\vee}^r (b {\vee}^r c) &= a {\vee}^r (\phi^{-1}(\phi(b) \vee \phi(c))	\\
	                                                &= \phi^{-1} (\phi(a) \vee \phi(\phi^{-1}(\phi(b) \vee \phi(c))))	\\
	                                                &= \phi^{-1} (\phi(a) \vee (\phi(b) \vee \phi(c)))	\\
	                                                & = \phi^{-1} ((\phi(a) \vee \phi(b)) \vee \phi(c))	\\
	                                                &= \phi^{-1} (\phi(\phi^{-1}(\phi(a) \vee \phi(b))) \vee \phi(c))	\\
	                                                &= \phi^{-1} (\phi(a {\vee}^r b)) \vee \phi (c)) = (a {\vee}^r b ) {\vee}^r c.
\end{aligned}$

$a {\wedge}^r (b {\wedge}^r c) = (a {\wedge}^r b) {\wedge}^r c$ holds in the same way. \\

\noindent  3. 
$\begin{aligned}[t]
	a {\wedge}^r (b {\vee}^r c) &= a {\wedge}^r \phi^{-1}(\phi(b) \vee \phi(c))	\\
	                                                      &=\phi^{-1}(\phi(a) \wedge \phi(\phi^{-1}(\phi(b) \vee \phi(c))))	\\
	                                                      &=\phi^{-1}(\phi(a) \wedge(\phi(b) \vee \phi(c)))	\\
	                                                      &=\phi^{-1}( (\phi(a) \wedge \phi(b)) \vee (\phi(a) \wedge \phi(c))).
\end{aligned}$

$\begin{aligned}
	(a {\wedge}^r b) {\vee}^r (a {\wedge}^r c) &= (\phi^{-1}(\phi(a) \wedge \phi(b))) {\vee}^r (\phi^{-1}(\phi(a) \wedge \phi(c)))	\\
	&= \phi^{-1}(\phi(\phi^{-1}(\phi(a) \wedge \phi(b))) \vee \phi(\phi^{-1}(\phi(a) \wedge \phi(c))))	\\
	&= \phi^{-1} (\phi(a) \wedge \phi(b)) \vee (\phi(a) \wedge \phi(c)).
\end{aligned}$

Hence \;$a {\wedge}^r (b {\vee}^r c) = (a {\wedge}^r b) {\vee}^r (a {\wedge}^r c)$ holds. 

$a {\vee}^r (b {\wedge}^r c) = (a {\vee}^r b) {\wedge}^r (a {\vee}^r c)$ holds in the same way. \\

\noindent 4. Let $I^r = \phi^{-1}(I^c)$, $O^r = \phi^{-1}(O^c)$. 
Then,

$\begin{aligned}
	a {\vee}^r O^r &= \phi^{-1}(\phi(a) \vee \phi(\phi^{-1}(O^c)))	\\
	                            &=\phi^{-1}(\phi(a) \vee O^c) = \phi^{-1}(\phi(a))=a
\end{aligned}$

$\begin{aligned}
	a {\wedge}^r I^r &= \phi^{-1}(\phi(a) \wedge \phi(\phi^{-1}(I^c)))	\\
	&=\phi^{-1}(\phi(a) \wedge I^c) = \phi^{-1}(\phi(a))=a, \text{\; for } {}^{\forall}a \in P(S).
\end{aligned}$

\noindent 5. For ${}^{\forall}a \in P(S)$, let $a'= \phi^{-1}((\phi(a))')$. Then,

$\begin{aligned}
	a {\vee}^r a' &=\phi^{-1}(\phi(a) \vee \phi(\phi^{-1}(\phi(a))'))	\\
	&=\phi^{-1}(\phi(a) \vee (\phi(a))') = \phi^{-1}(I^c)=I^r,
\end{aligned}$

$\begin{aligned}
	a {\wedge}^r a' &= \phi^{-1}(\phi(a) \wedge \phi(\phi^{-1}(\phi(a))'))	\\
	&= \phi^{-1}(\phi(a) \wedge (\phi(a))') = \phi^{-1}(O^c) =O^r.
\end{aligned}$

\noindent Hence, $P(S)$ is a Boolean algebra. $\square$ \\

\noindent We note that $O^r = \emptyset \in P(S)$ because $O^r = \phi^{-1}(O^c)= \phi^{-1}(\emptyset)= \emptyset$. 
In the above setting, we have the following proposition.

\begin{prop}
	For ${}^{\forall}a, b \in P(S)$, $\phi(a {\vee}^r b) = \phi(a) \vee \phi(b)$, $\phi(a {\wedge}^r b)= \phi(a) \wedge \phi(b)$, and $\phi(a')= \left( \phi(a) \right)'$.
	\label{prop-iso}
\end{prop}

\noindent Proof. 
$\phi(a {\vee}^r b)= \phi(\phi^{-1}(\phi(a) \vee \phi(b)))=\phi(a) \vee \phi(b)$. 
$\phi(a {\wedge}^r b)= \phi(\phi^{-1}(\phi(a) \wedge \phi(b)))=\phi(a) \wedge \phi(b)$. 
$\phi(a')=\phi(\phi^{-1}((\phi(a))')=(\phi(a))'$. $\square$

\vspace{3mm}

\noindent Hence, we have the following.

\begin{prop}
	The Boolean map $\phi: P(C) \to P(S)$ is an isomorphism of Boolean algebra. 
\end{prop}

\noindent Proof. $\phi: P(C) \to P(S)$ is bijective by Proposition \ref{prop-bi} and satisfies $\phi(a {\vee}^r b) = \phi(a) \vee \phi(b)$, $\phi(a {\wedge}^r b)= \phi(a) \wedge \phi(b)$ by Proposition \ref{prop-iso}. $\square$

\vspace{3mm}

\begin{figure}[h]
	\centering
	\includegraphics[width=100mm]{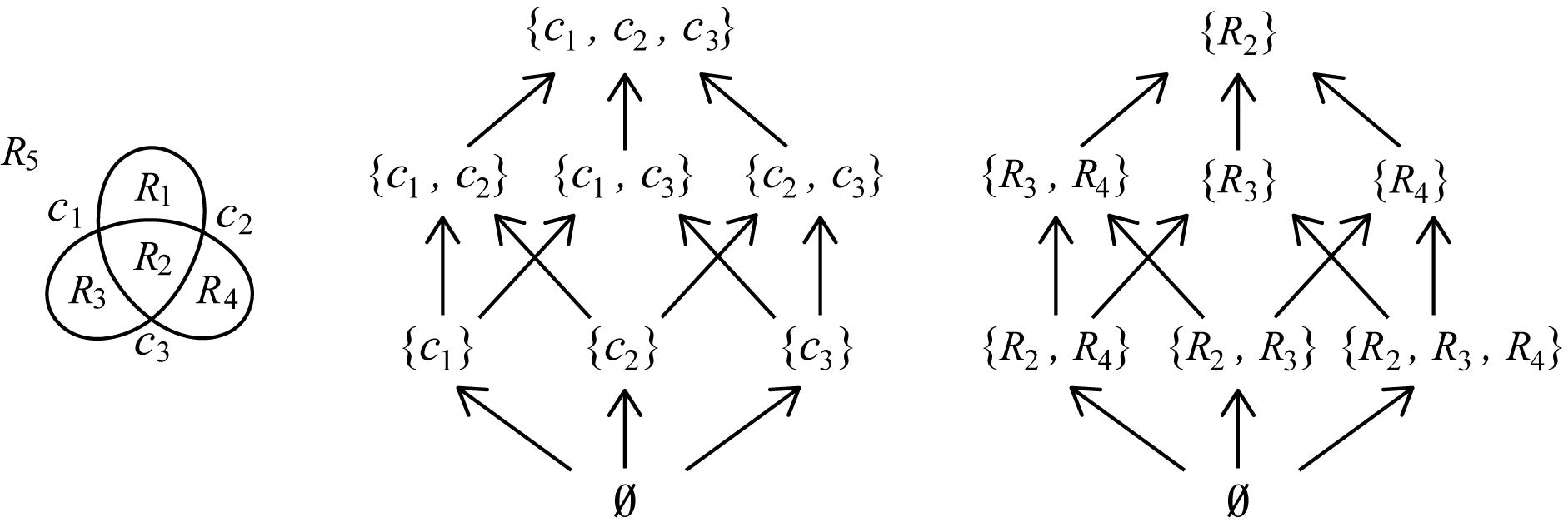}
	\caption{Boolean algebras $P(C)$ and $P(S)$ for $S= \{ R_2 , R_3 , R_4 \}$.}
	\label{fig-eq}
\end{figure}

\noindent The partial order ${\preceq}^r$ on $P(S)$ can be seen by Fact \ref{fact-b}, and we have the following correspondence. (See Figure \ref{fig-eq}.)

\begin{prop}
	$a {\preceq}^r b \Leftrightarrow \phi(a) \preceq \phi(b)$ for ${}^{\forall}a, b \in P(S)$. 
\end{prop}

\noindent Proof. 
$a {\preceq}^r b \Leftrightarrow a {\wedge}^r b' =O^r = \emptyset \Leftrightarrow \phi(a {\wedge}^r b')=\phi(\emptyset ) \Leftrightarrow \phi(a) \wedge \phi(b') = \emptyset \Leftrightarrow \phi(a) \wedge (\phi(b))' = O^c \Leftrightarrow \phi(a) \preceq \phi(b)$. $\square$

\section*{Acknowledgment}
	The second author's work was partially supported by JSPS KAKENHI Grant Number JP21K03263.

\end{document}